\documentclass[12pt]{amsart}

\usepackage[all]{xy}
\usepackage{amsmath, amsfonts, amssymb}

\newtheorem{thm}{Theorem}[section]
\newtheorem{lemma}[thm]{Lemma}
\newtheorem{pro}[thm]{Proposition}
\newtheorem{cor}[thm]{Corollary}
\newtheorem{defn}{Definition}[section]
\theoremstyle{remark}
\newtheorem{rem}[defn]{Remark}

\newtheorem{notation}[defn]{Notation}

\renewcommand{\Im}{\mathop{\mathrm{Im}}}
\newcommand{\pc}{\ensuremath{\mathbb{Z}/p\mathbb{Z}}}

\newcommand{\OJ}{\ensuremath{\Omega J_{p-1}(S^{2n})}}
\newcommand{\s}[1]{\ensuremath{S^{#1}}}

\newcommand{\ad}{\ensuremath{\mathrm{ad}}}
\newcommand{\lra}{\longrightarrow}
\newcommand{\dr}[3]{\ensuremath{#1\stackrel{#2}
{\longrightarrow}#3}}
\newcommand{\ddr}[5]{\ensuremath{#1\stackrel{#2}
{\longrightarrow}#3\stackrel{#4}{\longrightarrow}#5}}
\newcommand{\dddr}[7]{\ensuremath{#1\stackrel{#2}
{\longrightarrow}#3\stackrel{#4}{\longrightarrow}#5
\stackrel{#6}{\longrightarrow}#7}}
\newcommand{\ddddr}[9]{\ensuremath{#1\stackrel{#2}
{\longrightarrow}#3\stackrel{#4}{\longrightarrow}#5
\stackrel{#6}{\longrightarrow}#7}\stackrel{#8}{\longrightarrow}#9}

\newcommand{\hahc}{homotopy associative, homotopy commutative }
\newcommand{\OJp}{$\Omega J_{p-1}(S^{2n})$ }
\newcommand{\dls}{$\Omega^{2}S^{2n+1}$ }
\renewcommand{\L}{\ensuremath{\mathcal{L}}}
\newcommand{\U}{\ensuremath{\mathcal{U}}}

\begin{document}

\title[Universal H-spaces and the EHP spectral sequence]
{Universal homotopy associative, homotopy commutative $H$-spaces
and the $EHP$ spectral sequence}

\author{Jelena Grbi\' c}

\maketitle

\begin{abstract}
Assume that all spaces and maps are localised at a fixed prime
$p$. We study the possibility of generating a universal space
$U(X)$ from a space $X$ which is universal in the category of
homotopy associative, homotopy commutative $H$-spaces in the sense
that any map $f\colon X\lra Y$ to a homotopy associative, homotopy
commutative $H$-space extends to a uniquely determined $H$-map
$\overline{f}\colon U(X)\lra Y$. Developing a method for
recognising certain universal spaces, we show the existence of the
universal space $F_2(n)$ of a certain three-cell complex $L$.
Using this specific example, we derive some consequences for the
calculation of the unstable homotopy groups of spheres, namely, we
obtain a formula for the $d_1$-differential of the $EHP$-spectral
sequence valid in a certain range.
\end{abstract}

\section{Introduction}

In this paper we investigate the possibility of generating an
$H$-space out of its subspace and derive some consequences for the
calculations of the unstable homotopy groups of spheres.
Particular emphasis is laid on finding a formula for the first
differential $d_{1}$ in the $EHP$ spectral sequence which is valid
in a certain range.

Universal objects have always been of great interest to
mathematicians in  different mathematical disciplines.
Specifically, in homotopy theory one of the first universal spaces
is given in the category of homotopy associative $H$-spaces by the
James construction. Namely, if $X$ is a topological space with a
non-degenerate basepoint $e$ and a compactly generated topology,
then the James construction $J(X)$ on $X$ is the free topological
monoid generated by $X$ subject to the single relation that the
basepoint $e$ is the unit. The James construction satisfies the
following universal property: if $Y$ is a homotopy associative
$H$-space and a map $f\colon X\lra Y$ is given, then $f$ extends
to an $H$-map $\tilde{f}\colon J(X)\lra Y$ which is unique up to
homotopy.

Our aim is to find further information concerning universal
spaces. Indeed, one can ask what happens when the category of
homotopy associative $H$-spaces is replaced with the category of
\hahc ${H\text{-spaces}}$. Do universal spaces exist in this
category and if they do, is there a general construction of them?
While we do not answer the above questions in full generality, we
do present a new method for constructing universal spaces, and
describe some new universal spaces and their applications in
homotopy theory. For the sake of completeness we include all the
details as well as a brief overview of known examples of universal
spaces.

The method and its applications are $p$-local, so throughout the
remainder of the paper all spaces and maps have been localised at
a prime $p$. Also we assume that all spaces are pointed, connected
topological spaces with the homotopy types of finite type
$CW$-complexes. Unless otherwise indicated, the ring of homology
coefficients will be $\pc$ and $H_*(X;\pc )$ will be written as
$H_*(X)$.

In particular, aiming for a better understanding of the double
loop space on an odd sphere, we consider the first few spaces
appearing in Selick's filtration \cite{Se} of~$\Omega^2 S^{2n+1}$.
The main object of concern is a space $F_2(n)$ which is defined as
the homotopy pullback
\begin{equation}
\label{F_2(n)}
\xymatrix{F_{2}(n)\ar[d]\ar[r] & \Omega^{2} S^{2n+1}\ar[d]^{\Omega H}\\
S^{2np-1}\ar[r]^{E^2} & \Omega^{2} S^{2np+1},}
\end{equation}
where $H\colon\Omega S^{2n+1}\lra \Omega S^{2np+1}$ is the
James-Hopf invariant map and $E^2$ is the double suspension.

Let $L$ denote the $(2np-1)$-skeleton of $F_2(n)$. The central
result of the paper is as follows.

\begin*{thm}[A] The space $F_2(n)$ is the universal space of $L$
in the category of homotopy associative, homotopy commutative
$H$-spaces.
\end*{thm}
By $F_2(n)$ being a universal space in the category of \hahc
$H$-spaces, we mean that $F_2(n)$ has the two properties:
\begin{itemize} \item [(i)]\label{hahc} the space $F_2(n)$ is a
homotopy associative and homotopy commutative ${H\text{-space}}$;
\item [(ii)] \label{unispL} let $Y$ be a homotopy associative,
homotopy commutative $H$-space and $f\colon L\lra Y$ a given map.
Then there exists an $H$-map $\overline{f}\colon F_{2}(n)\lra Y$
extending $f$, and it is a unique up to homotopy,.
\end{itemize}
The universal property of $F_2(n)$ can be used to obtain useful
information about the $d_{1}$-differential in the $EHP$ spectral
sequence (cf. \cite{Ra:cc}). This spectral sequence is used to
calculate the homotopy groups of spheres (cf. \cite{To}), and is
made up by interlocking the classical $EHP$ sequences for various
spheres together. These $EHP$ sequences are long exact sequences
of the homotopy groups of spheres induced by the two classical
homotopy fibrations:
\[
\Omega^{2}S^{2np+1} \stackrel{P}{\longrightarrow}
J_{p-1}(S^{2n})\stackrel{E}{\longrightarrow}\Omega
S^{2n+1}\stackrel{H}{\longrightarrow}\Omega S^{2np+1}
\]
and
\[
\Omega^{2}S^{2np-1}\stackrel{P}{\longrightarrow}
 S^{2n-1}\stackrel{E}{\longrightarrow}\Omega
 J_{p-1}(S^{2n})\stackrel{H}{\longrightarrow}\Omega S^{2np-1},
\]
where $J_k(X)$ is the $k^{th}$-stage of the James construction on
a topological space $X$, and the map $H\colon\dr{\Omega
J_{p-1}(S^{2n})}{}{\Omega S^{2np-1}}$ is the Toda-Hopf invariant
map. In the literature the Toda-Hopf map is sometimes denoted by
$T$ to distinguish it more explicitly from the James-Hopf
invariant map.

Namely, for $p>2$ there are spectral sequences with $E^{k,
2n+1}_1=\pi_{k+2n+1}(S^{2np+1})$ and
$E^{k,2n}_1=\pi_{k+2n}(S^{2np-1})$. The spectral sequence with
$E^{k,n}=0$ for $n>j$ converges to $\pi_*(S^j)$ if $j$ is odd and
to $\pi_*(J_{p-1}(S^j))$ if $j$ is even (cf.\cite{Ra:cc}).

A key calculation in the $EHP$ spectral sequence is that of the
first differential $d_{1}\colon\pi_{r+2}(S^{2np+1})\lra
\pi_{r}(S^{2np-1})$. The ${d_1\text{-differential}}$ is induced by
the composition $\Omega^{3}S^{2np+1}\stackrel{\Omega
P}{\longrightarrow}\Omega
J_{p-1}(S^{2n})\stackrel{H}{\longrightarrow}\Omega S^{2np-1}$. In
the metastable range, the ${d_{1}\text{-differential}}$ is
completely determined by Toda's formula \cite{To}:
$d_{1}(E^{2}x)=p\cdot x$. ${\text{In \cite{Gr1}}}$ Gray
constructed a map $\varphi_{n}\colon\Omega^{2} S^{2np+1}\lra
S^{2np-1}$ with the property that ${\Omega\varphi_{n}
=H\circ\Omega P}$. Using the existence of the map  $\varphi_n$,
Gray \cite{Gr2} developed a formula for the
${d_{1}\text{-differential}}$ applicable to some elements which
are not double suspensions. In~\cite{Gr2} Gray showed that $\Omega
J_{p-1}(S^{2n})$ for $p>2$ is universal in the category of \hahc
$H$-spaces, with its generating subspace being the
$(2np-2)$-skeleton. He used this to extend the formula involving
$d_1$-differential. We shall extend further Gray's result using
the universality of $F_2(n)$.
\begin{notation}
In the case of odd spheres localised at a prime $p\geq 3$, the
degree $p$ map and the $p$-power map are homotopic and are
commonly denoted by $p\colon\dr{S^{2n+1}}{}{S^{2n+1}}$. Denote by
$S^{2n+1}\{p\}$ the homotopy fibre of the degree $p$ map on
$S^{2n+1}$. The homotopy cofibre of the degree $p^r$ map on
$S^{n-1}$ for $n\ge 2$ and $ r\ge 1$ is called the \emph{
$n$-dimensional $p$-primary Moore space} and is denoted by
$P^n(p^r)$. For any topological space $X$, the adjoint of the
identity map on $\Sigma X$ gives a map $E\colon X\lra \Omega\Sigma
X$, usually called the \emph{suspension map}. The \emph{double
suspension} $E^2\colon S^{2n-1}\lra\Omega^2 S^{2n+1}$ is the
double adjoint of the identity map on $S^{2n+1}$. Its homotopy
fibre is denoted by $W_n$. In \cite{Gr1} Gray has shown that $W_n$
has a classifying space $BW_n$.
\end{notation}
Using the definition of $F_2(n)$ as the homotopy pullback in diagram
\eqref{F_2(n)}, there exist two fibration sequences analogous to the
classical $EHP$ fibrations.
\begin*{pro}[B] There exist two homotopy
fibration sequences:
\[
 \dddr{W_{np}}{P}{F_2(n)}{E}{\Omega^2
S^{2n+1}}{H}{BW_{np}}
\]
and
\[
\dddr{\Omega S^{2np-1}\{
p\}}{P}{S^{2n-1}}{E}{F_2(n)}{H}{S^{2np-1}\{p\}}.
\]
\end*{pro}
Weaving together the long exact sequences of homotopy groups
induced by the above fibration sequences, we obtain an analogous
$EHP$ spectral sequence to the classical one.

For $p>2$ there are spectral sequences with $E^{k,
2n+1}_1=\pi_{k+2n+1}(W_{np})$ and
$E^{k,2n}_1=\pi_{k+2n}(S^{2np-1}\{ p\})$. The spectral sequence
with $E^{k,n}=0$ for $n>j$ converges to $\pi_*(S^j)$ if $j$ is odd
and to $\pi_*(F_2(\frac{j}{2}))$ if $j$ is even.

Using the universality of $F_{2}(n)$, we obtain the following
formula involving the $d_1$-differential.
\begin*{thm}[C] \label{formulad{1}} The composition
\[
F_2(np)\stackrel{E}{\longrightarrow}
\Omega^2S^{2np+1}\stackrel{\varphi_n}{\lra}S^{2np-1}\stackrel{E}{\lra}
 F_2(np)
\]
is the $p$-power map if either
\begin{enumerate}
 \item [(i)] $x\in\pi_*(F_2(np))$ is an element which is a lift through ${H\colon
 F_2(np)\lra S^{2np^2-1}\{p\}}$ of an element in the image of ${E\colon P^{2np^2-2}(p)\lra \Omega
 S^{2np^2-1}\{p\}}$
\end{enumerate}
or
\begin{enumerate}
 \item [(ii)] the composition $E\circ\varphi_n\circ E$ is an $H$-map.
\end{enumerate}
\end*{thm}

\begin*{cor}[D] Passing to homotopy groups under
either assumption (i) or assumption (ii) of Theorem~C, the formula
for the $d_1$-differential takes the form
\begin{equation}
\label{formularem} Ed_1(Ex)=p\cdot x,
\end{equation}
extending Gray's analogous result \cite{Gr2} with respect to
$\OJ$.
\end*{cor}

Theorem~C (ii) brings up once more the question whether
$\varphi_n$ is an ${H\text{-map}}$. This has not been proved as
yet. As a first step in that direction we have the following
result.

\begin*{pro}[E] \label{restphi} Restricted to the
$(2np^3-4)$-skeleton the composite
\[
F_2(np)\lra\Omega^2 S^{2np+1}\stackrel{\varphi_n}{\lra}
S^{2np-1}\stackrel{E}{\lra}F_2(np)
\]
is homotopic to the $p$-power map restricted to
$(F_2(np))_{(2np^3-4)}$.
\end*{pro}
Theorem~C and Proposition~E imply that formula~\eqref{formularem}
for the $d_1$-differential is valid up to the $(2np^3-4)$-skeleton
of $F_2(np)$. Since the $(2np^3-4)$-skeleton of $F_2(np)$ is
homotopy equivalent to the $(2np^3-4)$-skeleton of $\Omega^2
S^{2np+1}$, formula~\eqref{formularem} determines the composite of
the $d_1$-differential with $E\colon S^{2np-1}\lra F_2(np)$
restricted to the $(2np^3-4)$-skeleton of $\Omega^2S^{2np+1}$.

This improves Gray's result in the sense that
formula~\eqref{formularem} is not anymore valid just for
$(2np^2-2)$-skeleton but for $(2np^3-4)$-skeleton of
$\Omega^2S^{2np+1}$.
\section{Universal homotopy associative, homotopy commutative
 $H$-spaces}
\subsection{Preliminary definitions and known results}
This section deals with the concept of universal spaces, sometimes
called generated spaces. It sets up a new approach for
constructing universal spaces in the category of homotopy
associative, homotopy commutative $H$-spaces.
\begin{defn}
\label{generate} Let $X$ be a space, $G(X)$ an $H$-space and
${f\colon X\lra G(X)}$ a continuous map. $(X,f)$ will be called a
\textit{generator} for $G(X)$ if and only if the following
universal property is satisfied: for each $H$-space $Y$, the
mapping
\[
 \dr{\mbox{\{homotopy classes of $H$ maps $G(X)\rightarrow Y$\}}}
 {f^{\ast}}{[X,Y]}
\]
 is a bijection.
\end{defn}
\begin{rem} A space $G(X)$ satisfying the stated universal property is
called the {\it universal space of $X$}. Looking at the universal
property it is easy to notice that $G(X)$ is unique.
\end{rem}
One of the main examples of generated spaces is provided by the
James construction, which we described in the Introduction. In the
sense of Definition~\ref{generate}, the James construction has the
following universal property: if $Y$ is any homotopy associative
$H$-space and $f\colon X\lra Y$ is any continuous map, then $f$
has an extension to a unique $H$-map $\tilde{f}\colon J(X)\lra Y$.
We call $\tilde{f}$ the \textit{multiplicative extension of} $f$.
In particular, if $X$ is path connected, then the suspension map
$E\colon X\lra \Omega\Sigma X$ extends to a multiplicative
homotopy equivalence $\widetilde{E}\colon J(X)\lra\Omega\Sigma X$.
Therefore, the James construction on $X$ can be regarded as a
combinatorial model for $\Omega\Sigma X$. In this paper we
identify $J(X)$ with $\Omega\Sigma X$.


We can slightly modify the definition of the generating complex
asking for both spaces $G(X)$ and $Y$ to belong to the category of
\hahc $H$-spaces. With this additional condition Definition
\ref{generate} can be reformulated as follows.

\begin{defn}
\label{defuni} A \textit{universal space $U(X)$} of a space $X$ is
a \hahc $H$-space (localised at $p$) together with a map $i\colon
X\lra U(X)$ such that the following \textit{universal property}
holds:
\begin{itemize}
\item [] if $Y$ is a \hahc ${H\text{-space}}$ (localised at $p$)
and ${f\colon X\lra Y}$ is any map, then $f$ extends to a unique
${H\text{-map}}$ ${\overline{f}\colon U(X)\lra Y}$.
\end{itemize}
\end{defn}
Observe that the mapping induced by the definition of universal
spaces
\[
 \dr{\mbox{\{homotopy classes of $H$ maps $U(X)\rightarrow Y$\}}}
 {}{[X,Y]}
\]
 is a bijection.
\begin{rem}
Notice that the universal space $U(X)$ in the category of homotopy
associative, homotopy commutative $H$ spaces is defined by a space
localised at a prime $p$ rather then by an integral,
non-localised, space. The only reason for considering a localised
version of universal spaces is the existences of a better
developed machinery that can be applied to $p$-localised spaces in
order to describe their properties.
\end{rem}

\subsubsection*{Known examples of universal spaces in the
category of \hahc $H$-spaces.}
\addcontentsline{toc}{subsubsection}{Known examples of universal
spaces in the category of \hahc $H$-spaces.}

\noindent{\it1.} The simplest example of universal spaces appears
when we look for a universal space of a torsion free space, that
is, of a rational space. The universal space of a rational space
$X$ always exists and is given by the infinite symmetric product
$SP^\infty(X)$. This fact was proved by Dold and Thom \cite{DT},
although they did not use the terminology of universal spaces.

\noindent{\it 2.} The fact that some examples of $U(X)$ exist,
even with the additional homology condition, can easily be seen in
the case of spheres. The universal space $U(S^{2n})$ of an even
dimensional sphere for $p\geq 3$ is already given by the James
construction on $S^{2n}$ since $\Omega S^{2n+1}$ is homotopy
commutative as well as homotopy associative \cite{Ad:sp}. The
universal space $U(S^{2n-1})$ of an odd dimensional sphere is
almost implicitly given by Serre's splitting \cite{Serre} at odd
primes of the loop space on an even dimensional sphere $ \Omega
S^{2n}\stackrel{\simeq}{\longrightarrow} S^{2n-1}\times \Omega
S^{4n-1} $. The equivalence is obtained via the product of the map
$E \colon \s{2n-1}\lra \Omega \s{2n}$ with the loop map on the
Whitehead product $\omega\colon \s{4n-1}\lra\s{2n}$ of the
identity map on $S^{2n}$ with itself. For $p\geq 5$, an odd
dimensional sphere $S^{2n-1}$ is homotopy associative, homotopy
commutative, and $U(S^{2n-1})=S^{2n-1}$. A more detailed proof for
the universal spaces of spheres reader can find in \cite{G}.

\noindent{\it 3.} The next example of universal spaces occurs in
the work of Cohen, Moore and Neisendorfer \cite{CMNt} on splitting
the loop spaces on even dimensional odd primary Moore spaces
${\Omega P^{2n+2}(p^r)\stackrel{\simeq}{\lra} S^{2n+1}\{ p^r\}
\times \Omega (\vee_\alpha P^{n_\alpha}(p^r))}$. In this case,
Gray showed \cite{Gr2} that $U(P^{2n+1}(p^r))=S^{2n+1}\{ p^r\}$,
the homotopy fibre of the degree $p^r$ map on $S^{2n+1}$.

\noindent{\it 4.} The most difficult known examples of $U(X)$ in
the sense that they do not retract off the James construction
$J(X)$ on $X$ are those $H$-spaces whose homology rings are the
abelianisation of the homology ring of the loop space on odd
dimensional Moore spaces at odd prime. These spaces have been
constructed by Anick and Gray \cite{Ad}, \cite{AG} and their
universality has been shown by Theriault \cite{Th2}. In this case
the universal space of $P^{2n}(p^r)$ does not split off the James
construction on $P^{2n}(p^r)$ unlike the other cases above.

\noindent{\it 5.} An example of a universal space that is closely
related to the universal space we are interested in considers a
loop space $\Omega J_{p-1}(S^{2n})$. Let $K_{n}$ be the
$(2np-2)$-skeleton of \OJp. Gray \cite{Gr2} has shown that $\Omega
J_{p-1}(S^{2n})$ is homotopy commutative as well as homotopy
associative, due to the loop space structure, and that it is the
universal space of $K_n$.

\subsection{A homological filtration of the double suspension of
odd dimensional spheres} In order to construct a spectral sequence
concerning the double suspension Selick \cite{Se} looked at a
certain filtration of $\Omega^{2}S^{2n+1}$. We recall it here
briefly. Looping the filtration
\[
\{pt\}\longrightarrow J_{1}(S^{2n})\longrightarrow
\ldots\longrightarrow
J_{k}(S^{2n})\longrightarrow\ldots\longrightarrow
J(S^{2n})\simeq\Omega\Sigma S^{2n}
\]
of the James construction on $S^{2n}$ gives a filtration
\[
\{pt\}\longrightarrow \Omega J_{1}(S^{2n})\longrightarrow
\ldots\longrightarrow \Omega
J_{k}(S^{2n})\longrightarrow\ldots\longrightarrow
\Omega^{2}S^{2n+1}
\]
of $\Omega^{2}S^{2n+1}$.

It is well known (cf. \cite{CML}) that
\[
H_{\ast}(\Omega^{2}S^{2n+1})\cong\bigotimes^{\infty}_{j=0}\Lambda
(x_{2np^{j}-1})\otimes\bigotimes^{\infty}_{j=1}
\mathbb{Z}/p\mathbb{Z} [y_{2np^{j}-2}]
\]
with $\beta x_{2np^{j}-1}=y_{2np^{j}-2}$. Applying the Serre or
Eilenberg-Moore spectral sequence to the homotopy fibration
sequence
\[
\Omega J_{p^k-1}(S^{2n})\lra \Omega^2 S^{2n+1}\lra \Omega^2
S^{2np^k+1},
\]
it follows that
\[
H_{\ast}(\Omega
J_{p^{k}-1}(S^{2n}))\cong\bigotimes^{k-1}_{j=0}\Lambda
(x_{2np^{j}-1})\otimes\bigotimes^{k}_{j=1} \mathbb{Z}/p\mathbb{Z}
[y_{2np^{j}-2}].
\]
Therefore, the filtration
\begin{equation}
\label{filtration}
 \{pt\}\lra
\Omega J_{p-1}(S^{2n})\longrightarrow \ldots\longrightarrow \Omega
J_{p^{k}-1}(S^{2n})\longrightarrow\ldots\longrightarrow
\Omega^{2}S^{2n+1}
\end{equation}
gives a natural filtration of $H_{\ast}(\Omega^{2}S^{2n+1})$. The
aim is to obtain a full homological filtration of
$\Omega^2S^{2n+1}$. To refine filtration (\ref{filtration}) we
shall insert spaces, denoted by $F_{2k}(n)$, between $\Omega
J_{p^{k}-1}(S^{2n})$ and $\Omega J_{p^{k+1}-1}(S^{2n})$ requiring
that their homology is
\[
H_{\ast}(F_{2k}(n))\cong\bigotimes^{k}_{j=0}\Lambda
(x_{2np^{j}-1})\otimes\bigotimes^{k}_{j=1} \mathbb{Z}/p\mathbb{Z}
[y_{2np^{j}-2}].
\]
Selick's filtration
\begin{equation}
\label{filtration1} F_{-1}(n)\longrightarrow
F_{0}(n)\longrightarrow F_{1}(n)\rightarrow\ldots\rightarrow
F_{k}(n)\longrightarrow\ldots\longrightarrow F_{\infty}(n)
\end{equation}
is defined inductively on $n$ and $k$. Set the initial data in a
way that $F_{-1}(n)=\{pt\}$, $F_{0}(n)=S^{2n-1}$ and
$F_{\infty}(n)=\mbox{\dls}$.  Assume by induction on $k$ that for
all $q\le k$ and for all $n$ the spaces $F_{q}(n)$ have been
constructed along with maps
\[
F_{-1}(n)\lra F_{0}(n)\longrightarrow
F_{1}(n)\longrightarrow\ldots\longrightarrow
F_{k}(n)\longrightarrow F_{\infty}(n).
\]
The space $F_{k+1}(n)$ is defined as the homotopy pullback
\[
\xymatrix{ F_{k+1}(n) \ar[r] \ar[d] & F_{\infty}(n) \ar[d]^{\Omega
H}\\
F_{k-1}(np)\ar[r] & F_{\infty}(np) ,}
\]
where $H\colon\Omega S^{2n+1}\longrightarrow \Omega S^{2np+1}$ is
the James-Hopf invariant map. Being a homotopy pullback of
$H$-spaces and $H$-maps each $F_{k}(n)$ is an $H$-space and all
the induced maps in the pullback diagram are $H$-maps. Toda's
calculations \cite{Tod} show that $F_{2k-1}(n)\simeq \Omega
J_{p^{k}-1}(S^{2n})$. Therefore, we conclude that in filtration
(\ref{filtration1}) all the spaces are $H$-spaces and all the maps
are $H$-maps.

In particular, $F_{2}(n)$ is defined as the homotopy pullback
\[
\xymatrix{F_{2}(n)\ar[d]\ar[r] & \Omega^{2} S^{2n+1}\ar[d]^{\Omega H}\\
S^{2np-1}\ar[r]^-{E^2} & \Omega^{2} S^{2np+1}.}
\]
In this way, $F_{2}(n)$ can be considered as the third stage in an
approximation of $\Omega^2S^{2n+1}$.

Therefore, the knowledge of $F_{2}(n)$ can be used to better
understand the map $\varphi_n\colon\Omega^2 S^{2np+1} \lra
S^{2np-1}$, discussed in the Introduction, and to extend the
calculations of the $d_1$-differential in the $EHP$ spectral
sequence.

\subsection{The strategy for the proof of Theorems B and C}

In this short section we outline the strategy for the proof of
Theorems B and C in order to facilitate the reading and the
understanding of the main ideas. We briefly give the major facts
and procedures.

Consider the homology of $F_{2}(n)$, that is,
\begin{equation}
\label{fhomology}
 H_{\ast}(F_{2}(n))\cong \Lambda (x_{2n-1},
x_{2np-1})\otimes \mathbb{Z}/p\mathbb{Z} [y_{2np-2}],
\end{equation}
the symmetric algebra on three generators $x_{2n-1},y_{2np-2},
x_{2np-1}$.

Denote by $K$ the $(2np-2)$-skeleton of $\OJ$. Then $K$ maps to
$F_2(n)$ by the composition
\[
K\lra \OJ =F_1(n)\lra F_2(n)
\]
implying that the pair $(F_2(n), K)$ is $2np-2$ connected. Let
\[
\dr{(h,\ q)\colon (e^{2np-1}, \s{2np-2})}{}{(F_2(n), K)}
\]
represent a generator of
\[
\pi_{2np-1}(F_2(n), K)\cong H_{2np-1}(F_2(n), K;\ \mathbb{Z}).
\]
Define a space $L$ as the homotopy cofibre of $q$, namely, as the
three-cell complex $L=K\cup_q e^{2np-1}$. Hence $L$ is the
$(2np-1)$-skeleton of $F_2(n)$, where each cell corresponds to a
generator of $S(x_{2n-1}, y_{2np-2}, x_{2np-1})$.

The general idea in constructing the universal space of $L$ is to
look first at the James construction on $L$. Since $\Omega\Sigma
L$ is not homotopy commutative (if it were, its homology would be
a symmetric algebra) it cannot be the universal space of $L$.
Therefore, for a suitable candidate for the universal space we aim
for a homotopy associative, homotopy commutative retract of
$\Omega\Sigma L$ containing $L$ as a subspace. Since the universal
space of $L$ is by definition homotopy commutative, its homology
is commutative, so an initial guess for the homology ring would be
the abelianisation of $T(\widetilde{H}_*(L))$, that is the
symmetric algebra $S(\widetilde{H}_*(L))$. Our method for
constructing the universal space of $L$ involves three steps. In
the first two steps we are constructing a suitable retract of
$\Omega\Sigma L$. We start by obtaining a decomposition of
$H_*(\Omega\Sigma L)$, where one of the factors is a symmetric
algebra. To do this we regard $H_*(\Omega\Sigma L)$ as a universal
enveloping algebra and apply the theory of Lie algebras.
\begin{lemma}
 \label{homOSL} There is an isomorphism
\[
H_*(\Omega\Sigma L)\cong T(W)\otimes S(x_{2n-1}, y_{2np-2},
x_{2np-1}),
\]
of left $T(W)$-modules and right $S(x_{2n-1}, y_{2np-2},
x_{2np-1})$-comodules, where $W$ is a free $\pc$-module.
\end{lemma}
The construction is then based on a geometric realisation of the
homological decomposition, where each factor is realised as a
space of a certain principal fibration.
\begin{pro}
\label{fibre} There is a fibration sequence
\begin{equation}
\label{fib} \ddddr{\Omega\Sigma R}{\Omega\omega}{\Omega\Sigma
L}{}{F_2(n)}{}{\Sigma R}{\omega}{\Sigma L},
\end{equation}
where $\Sigma R$ is a wedge of spheres and Moore spaces, and
$\omega$ is a sum of Whitehead products and mod~$p$ Whitehead
products. Moreover, $\Omega\omega$ has a left homotopy inverse.
\end{pro}
\begin{cor}
There is a homotopy decomposition
\[
\Omega\Sigma L\simeq F_2(n)\times \Omega\Sigma R.
\]
\end{cor}
If Proposition \ref{fibre} is granted, then \eqref{fhomology}
shows that $F_2(n)$ corresponds to a geometric realisation of the
symmetric algebra $S(x_{2n-1}, y_{2np-2}, x_{2np-1})$. In the
third and last step we shall further examine fibration~\eqref{fib}
in order to prove that $F_2(n)$ is the universal space space of
$L$ in the category of \hahc $H$-spaces, that is, we shall prove
Theorems~B and~C.

\subsection{The construction of a map $\omega\colon\Sigma R\lra\Sigma L$}

We want to construct a map $\omega$ which takes a certain wedge of
spheres and Moore spaces $\Sigma R$ into $\Sigma L$ by means of
Whitehead products and mod $p$ Whitehead products in such a way
that the homology map $(\Omega\omega)_{\ast}$ is an isomorphism
onto the subalgebra generated by Lie brackets in
$H_{\ast}(\Omega\Sigma L)$. Therefore, $\omega$ realises a
decomposition of the homology of $\Omega\Sigma L$ given as a
product of $H_*(\Omega\Sigma R)$ and some symmetric algebra.

A purely algebraic approach gives prominent insight into this
problem. The following analysis has been drawn from the work of
Cohen, Moore and Neisendorfer~\cite{CMNt}.

\begin{notation}
Throughout this section, let $R$ be a commutative ring in which 2
is a unit. If $V$ is positively graded module over $R$, then
$\mathcal{L}\langle V\rangle$ is the free graded Lie algebra over
$R$ generated by $V$. In particular,
$\mathcal{L}=\mathcal{L}\langle x,u,v\rangle$ is the Lie algebra
over $R$ on three generators $x,u,v$ in degrees $2n-1$, $2np-2$,
and $2np-1$, respectively. Denote by $[\mathcal{L} ,\mathcal{L}]$
the commutator of $\mathcal{L}$. We regard $\L_{ab}\langle
V\rangle$ as the free graded abelian Lie algebra over $R$
generated by $V$. The universal enveloping algebra of
$\mathcal{L}$ is denoted by $U\mathcal{L}$. The universal
enveloping algebra $U\mathcal{L}\langle V\rangle$ of the free Lie
algebra $\mathcal{L}\langle V\rangle$ is isomorphic to the free
tensor algebra $T(V)$ generated by $V$, and hence $U\L_{ab}\langle
V\rangle$ is isomorphic to the free symmetric algebra $S(V)$
generated by $V$.
\end{notation}

By the definition of a free graded abelian Lie algebra there is a
short exact sequence of Lie algebras
\[
0\longrightarrow
[\mathcal{L},\mathcal{L}]\longrightarrow\mathcal{L}\longrightarrow
\L_{ab}\langle x,u,v \rangle\longrightarrow 0.
\]
According to \cite[Proposition 3.7]{CMNt}, this results in a short
exact sequence of homology Hopf algebras
\[
0\longrightarrow U[\mathcal{L},\mathcal{L}]\longrightarrow
U\mathcal{L}\longrightarrow U\L_{ab}\langle x,u,v
\rangle\longrightarrow 0,
\]
and an isomorphism
\begin{equation}
\label{Liedecomposition} U\mathcal{L}\cong U
[\mathcal{L},\mathcal{L}]\otimes U\L_{ab}\langle x,u,v \rangle
\end{equation}
of left $U[\mathcal{L}, \mathcal{L}]$-modules and right
$U\L_{ab}\langle x,u,v \rangle$-comodules. On the other hand,
whenever $\mathcal{L}$ is a connected free Lie algebra its
commutator $[\mathcal{L},\mathcal{L}]$, as a sub-Lie algebra of
$\mathcal{L}$, is a free Lie algebra $\mathcal{L}\langle W\rangle$
with $W$ a free $R$-module. Therefore it follows that
\[
U[\mathcal{L},\mathcal{L}]\cong T(W).
\]

\begin{proof}[Proof of Lemma~{\rm\ref{homOSL}}]
Notice that $H_*(L)$ is the free $\pc$--module with basis $x,u$
and $v$ in degrees $2n-1$, $2np-2$ and $2np-1$. Since $H_*(L)$ is
a trivial coalgebra, $H_*(\Omega\Sigma L)$ is the primitively
generated tensor Hopf algebra $T(x_{2n-1}, u_{2np-2}, v_{2np-1})$
generated by $x$, $u$ and $v$. Therefore $H_*(\Omega\Sigma L)$ can
also be considered as the universal enveloping Hopf algebra $U\L$
of the free Lie algebra $\L=\L\langle x, u, v\rangle$.

Applying decomposition \eqref{Liedecomposition} of the universal
enveloping algebra $H_*(\Omega\Sigma L)$, we get an isomorphism
\[
H_*(\Omega\Sigma L)\cong T(W)\otimes S(x_{2n-1}, u_{2np-2},
v_{2np-1}),
\]
of left $T(W)$-modules and  right $S(x_{2n-1}, u_{2np-2},
v_{2np-1})$-comodules, where $W$ is a free $\pc$-module. This
proves Lemma~\ref{homOSL}.
\end{proof}

We would like to have  better control over the free $\pc$--module
$W$. Recall that the Euler-Poincar\' e series of a free graded
module $V=\bigoplus_{n=0}^{\infty} V_n$ of finite type is the
power series
\[
\chi (V)=\sum^{\infty}_{n=0} (\mathop{\mathrm{rank}} V_{n})t^{n}.
\]

Suppose that $\L'\lra \L\lra\L''$ is a short exact sequence of
connected graded Lie algebras which are free finite type
$R$-modules. If $\L$ is a free Lie algebra $\L\langle V\rangle$,
then by \cite[Corollary 3.11]{CMNt} $\L'$ is a free Lie algebra
$\L\langle W\rangle$, where $W$ is a free finite type $R$-module
with the Euler-Poincar\'{e} series given by
\begin{equation}
\label{eulereq}
 \chi (W)=1-\chi (\U\L'')(1-\chi (V)).
\end{equation}

In our case $\mathcal{L}$ is the free Lie algebra on three
generators $x_{2n-1},$ $u_{2np-2}$ and $v_{2np-1}$, and we have
written the commutator $[\mathcal{L},\mathcal{L}]$ as the free Lie
algebra $\mathcal{L}(W)$ generated by a finite-type $\pc$-module
$W$.

\begin{lemma}
\label{euPO} The Euler-Poincar\'{e} series of
$[\mathcal{L},\mathcal{L}]$
has the following form\\
\hspace*{2cm}$ \chi(W)=
(t^{4n-2}+t^{2n(p+1)-3}+t^{2n(p+1)-2}+t^{2n(p+2)-3}+$\\
\hspace*{4cm}$t^{4np-3}+t^{4np-2}+t^{2n(2p+1)-4}+t^{2n(2p+1)-3})\sum^{\infty}_{k=0}
t^{2k(np-1)}.$
\end{lemma}
\begin{proof}
The Euler-Poincar\'{e} series for $W$ is given by
equation~\eqref{eulereq}, namely,
\begin{equation}
\label{Poincare} 1-\chi(W)=\chi(\mathcal{U}\L_{ab}\langle
x,u,v\rangle)(1-\chi(x,u,v)).
\end{equation}
Expanding the right side of equation (\ref{Poincare})gives
\[
1-\chi(W)=\frac{(1+t^{2n-1})(1+t^{2np-1})}{1-t^{2np-2}}(1-t^{2n-1}-t^{2np-2}-t^{2np-1}).
\]
Sorting out the previous expression the stated formula for $\chi
(W)$ is obtained.
\end{proof}
\begin{rem}
For $x,y\in\mathcal{L}$, the Lie bracket $[x,y]$ in $\mathcal{L}$
is also denoted by $\ad(x)(y)$. This notation can be generalised
by defining $\ad^0(x)(y)=y$ and for $k\geq 1$, inductively
defining $\ad^k(x)(y)=\ad(x)(\ad^{k-1}(x)(y))$ for each
$x,y\in\mathcal{L}$.
\end{rem}
For a Hopf algebra $A$, let $Q(A)$ denote the module of
indecomposables. Recall from \cite{CMNt} that whenever there is a
short exact sequence $\ddr{A'}{}{A}{}{A''}$ of homology Hopf
algebras the Lie bracket $[\ ,\ ]\colon A\otimes A\lra A$ induces
a left action of $A''$ on $Q(A')$.
\begin{lemma} [{\cite [Lemma 3.12]{CMNt}}]
 \label{action}
Let $A$ be a tensor algebra $T(V)$ and suppose that $A''$ is a
free module over a tensor subalgebra $T(U)$ of $A''$. Then $Q(A')$
is a free $T(U)$-module.
\end{lemma}
\begin{pro}
\label{basis} The commutator $[\L ,\L]$ of the free Lie algebra
$\L=\L\langle x,u,v\rangle$ is a free Lie algebra with the
following
generators\\
 $\ad^{k}(u)[x,x]$, $\ \ad^{k}(u)[x,u]$, $\
\ad^{k}(u)[x,v]$, $\ \ad^{k}(u)[u,v]$, $\ \ad^{k}(u)[v,v]$, $\
\ad^{k}(u)[x,[u,v]]$, $\ \ad^{k}(u)[x,[v,v]]$, $\
\ad^{k}(u)[v,[x,x]]$ for $k\geq 0$.
\end{pro}
\begin{proof}
The degrees of the rank 2 commutators are different, implying that
they are independent. For commutators of rank $3$, because of the
dimensional reason the only possible relation might hold between
$[u,[x,v]]$ and $[x,[u,v]]$. But by the Jacobi identity
$[u,[x,v]]= [[u,x],v]+[x,[u,v]]$ and the fact that the first of
these summands is non-zero, these two commutators are independent.
The higher rank commutators are obtained by iterating the
application of $\ad(u)$. Applying Lemma~\ref{action} to the tensor
subalgebra $T(u)\subset S(x,u,v)$, the action of $\ad(u)$ is free
on $W$, and therefore all the listed elements are independent. We
finish the proof by comparing the given set $W$ and the
Euler--Poincar\'{e} series given by Lemma \ref{euPO}.
\end{proof}

Now we bring topology into the picture. We are aiming to the
geometric realistion of the homology decomposition of
$\Omega\Sigma L$ given by Lemma \ref{homOSL}. For that we need to
recall some preliminary definitions.

Let $p$ be a prime and $n\geq 2$. The $n$-dimensional mod~$p$ {\it
Moore space} $P^n(p)$ is the cofibre of the degree $p$ map on
$S^{n-1}$. The {\it $n^{th}$ mod~$p$ homotopy set} of a pointed
space $X$ is the set of based homotopy classes
\[
\pi_n(X;\mathbb{Z}/p\mathbb{Z})=[P^n(p),X].
\]
The $n^{th}$ mod~$p$ homotopy set
$\pi_n(X;\mathbb{Z}/p\mathbb{Z})$ has a group structure for $n>2$
and an abelian group structure if $n>3$.

Let $G$ be a group-like space. The \textit{mod p Samelson product}
is a pairing
\[
\pi_k(G;\mathbb{Z}/p\mathbb{Z})\otimes\pi_l(G;\mathbb{Z}/p\mathbb{Z})
\lra \pi_{k+l}(G;\mathbb{Z}/p\mathbb{Z})
\]
whose value at $f\colon P^k(p)\lra G$ and $g\colon P^l(p)\lra G$
is given by the composite
\[
[f,g]\colon P^{k+l}(p)\lra\ddr{P^k(p)\wedge P^l(p)}{f\wedge
g}{G\wedge G}{[\,,\ ]}{G},
\]
 where $[\,,\ ]$ is a commutator in
$G$.

The Bockstein homomorphism $\beta\colon\pi_n(G; \pc)\lra
\pi_{n-1}(G;\pc)$ is induced by the composite
\[
\ddr{P^{n-1}(p)}{}{S^{n-1}}{}{P^n(p)}
\]
of the pinch map from with the inclusion. Because of the
dimensional reason the composition $\beta\circ\beta$ is trivial,
namely, the Bockstein homomorphism is a differential on the graded
module $\pi_*(G;\pc)$.

The mod~$p$ Samelson product is compatible with the Bockstein,
that is,
\[
 \beta [f,g]=[\beta f,g]+(-1)^k[f,\beta g].
\]
If $p\geq 5$, the mod~$p$ Samelson product and the Bockstein
operator give $\pi_*(G;\pc )$ a differential graded Lie algebra
structure \cite{Ne}.

The \textit{mod $p$ Whitehead product} is defined by taking the
adjoint of mod $p$ Samelson product.

There is a mod~$p$ Hurewicz homomorphism defined as follows. Let
$v$ be a generator of $H_n(P^n(p);\mathbb{Z}/p\mathbb{Z})$ and
$[l]\in\pi_n(X;\mathbb{Z}/p\mathbb{Z})$. The Hurewicz homomorphism
\[
\dr{h\colon\pi_n(X;\mathbb{Z}/p\mathbb{Z})}{}{H_n(X;\mathbb{Z}/p\mathbb{Z})}
\]
is defined by $h([l])=l_*(v)$.

Another important property of the $p$-primary Samelson product
describes its behavior with respect to the Hurewicz homomorphism.
As in an associative algebra, the Lie bracket in $H_*(G;\pc)$ is
defined by $[x,y]=xy-(-1)^{|x||y|}yx$ for $x,y\in H_*(G;\pc)$.
With $f$ and $g$ given as above, we have
\begin{equation}
\label{HurvsSamelson}
 h([f,g])=[h(f),h(g)].
\end{equation}
Hence for $p\geq 5$, the Hurewicz homomorphism is a Lie algebra
homomorphism from the $p$-primary homotopy groups of $G$ to the
Lie algebra of primitives in $H_*(G;\pc )$.

\begin{lemma}
\label{HurImage} Each of the basis elements in Proposition
\ref{basis} is in the image of the mod~$p$ Hurewicz homomorphism
in $U\L$ except $[x,x]$ which is spherical.
\end{lemma}
\begin{proof}
Let $\chi$ be the identity map on sphere $S^{2n-1}$, and $\L(\chi
)$ the Lie algebra generated by $\chi$. Then the Hurewicz
homomorphism $h$ maps $\chi$ onto $x$. The Hurewicz homomorphism
commutes with Samelson products and therefore maps $[\chi,\chi ]$
onto $[x,x]$. Now, let $\nu$ and $\mu$ denote the identity and the
Bockstein maps on $P^{2np-1}(p)$. Then $\beta (\nu)=\mu$ and
$\beta (\mu)=0$. Let $\mathcal{L}\langle\nu ,\mu\rangle$ be the
free differential graded Lie algebra over $\mathbb{Z}/p\mathbb{Z}$
generated by $\{ \nu ,\mu\}$ with $\beta$ as a differential and
the Lie structure given by the mod~$p$ Whitehead product. The
Hurewicz homomorphism acts as $h(\nu)=v$ and $h(\mu)=u$, where
$u,v\in U\mathcal{L}\cong H_*(\Omega\Sigma L)$. By Formula
\eqref{HurvsSamelson}, the mod~$p$ Hurewicz homomorphism commutes
with mod~$p$ Whitehead products. Now the Lemma can be proved by
induction on the length of the Lie bracket.
\end{proof}

To define the desired map $\omega$ take $R$ to be the wedge of
spheres and Moore spaces that correspond to the basis $W$ from
Proposition~\ref{basis} and map $\Sigma R$ by Whitehead products
and mod~$p$ Whitehead products to $\Sigma L$. More precisely, let
$W=\{y_{i}\}_{i\in I}\cup z$, $\alpha\colon S^{4n-2}\lra
\Omega\Sigma L $ the Samelson product whose Hurewicz image is $z$,
$\beta_i\colon P^{n_i}_{i}(p)\lra \Omega\Sigma L$ the mod~$p$
Samelson product whose Hurewicz image is $y_{i}$ and
$R=S^{4n-2}\bigvee\Bigl(\bigvee_{i\in I} P^{n_i}_{i}(p)\Bigr)$.
Define $\bar{\omega}\colon R\lra \Omega\Sigma L$ by taking the
coproduct of $\alpha$ and $\beta_i$ over $i\in I$. Now define the
map $\omega\colon\Sigma R\lra \Sigma L$ as the adjoint of
$\bar{\omega}$.
\begin{pro}
\label{omegaomega*} In $\pc$--homology, the induced map
\[
(\Omega\omega)_*\colon U [\L,\L]\lra U\L
\]
is a Hopf algebra monomorphism.
\end{pro}
\begin{proof}
Being an $H$-map into a homotopy associative $H$-space,
$\Omega\omega$ is determined by its restriction to $R$ by the
James construction. The restriction $\overline{\omega}\colon R\lra
\Omega\Sigma L$ is, by construction, the coproduct of Samelson and
mod $p$ Samelson products whose Hurewicz images biject onto the
basis $W$ of $[\L,\L]$. Extending multiplicatively
$\overline\omega$ to $\Omega\omega$, the homology map
$(\Omega\omega)_*$ sends $H_*(\Omega\Sigma R)$ isomorphically onto
the subalgebra $U[\L,\L]$ generated by Lie brackets in
$H_*(\Omega\Sigma L)$.
\end{proof}


\subsection{The homotopy fibre of $\omega$}

Denote by $F$ the homotopy fibre of $\omega$. In this section we
consider the fibration sequence
\begin{equation}
\label{fibre2} \ddddr{\Omega\Sigma R}{\Omega\omega}{\Omega\Sigma
L}{\partial}{F}{}{\Sigma R}{\omega}{\Sigma L}
\end{equation}
and show that $F$ is homotopy equivalent to $F_{2}(n)$. Firstly,
we show that $\Omega\omega$ has a left homotopy inverse. This
implies that there is a homotopy decomposition of $\Omega\Sigma L$
as $\Omega\Sigma L\simeq F\times\Omega\Sigma R $, and leads
towards the conclusion that $H_{\ast}(F)\cong S(x_{2n-1},
y_{2np-2}, x_{2np-1})$.

To begin, we need several of the properties of $L$ and
 $F_{2}(n)$. In the grand design we follow Selick's work~\cite{Se} on the
 double suspension. On our way through it, we prove new
 properties of these spaces which we shall require subsequently.

\begin{thm}[{[Selick~\cite{Se}]}]
\label{properties} The spaces $F_2(n)$ and $L$ have the following
properties:
\begin{enumerate}\renewcommand{\labelenumi}{(\roman{enumi})}
 \item As a Hopf algebra,
\[
H_{\ast}(F_{2k}(n))\cong \bigotimes^{k}_{j=0} \Lambda
(x_{2np^{j}-1})\otimes\bigotimes^{k}_{j=1} \mathbb{Z}/p\mathbb{Z}
[y_{2np^{j}-2}]
\]
with the generators primitive. The action of the Bockstein is
given by\\ ${\beta (x_{2np^{j}-1})=y_{2np^{j}-2}}$;
 \item If $n>1$, then $F_{k}(n)$ is atomic for $1\leq k\leq\infty$;
 \item Let $K$ be the bottom $2$-cell subcomplex of $L$. Then
 \[
  \Sigma K\wedge L\simeq \s{4n-1}\vee \s{2(np+n)-2}\vee
P^{2(np+n)-1}(p)\vee P^{4np-2}(p);
  \]
 \item \label{ssl} $\Sigma^{2} L\simeq S^{2n+1}\vee P^{2np+1}(p);$
 \item \label{sigmaf2} $\Sigma F_{2}(n)\simeq \Sigma L\vee M$, where $M$ is a wedge of Moore
spaces.
\end{enumerate}
\end{thm}

Particular emphasis should be placed upon the next property of the
space $L$ since it is going to play a crucial role in proving that
$F_{2}(n)$ is homotopy associative and homotopy commutative. While
the next result is probably known, to author's knowledge it does
not appear in the literature.
\begin{lemma}
\label{sll}The space $\Sigma L\wedge L$ is homotopy equivalent to
\[
S^{4n-1}\vee P^{2n(p+1)-1}(p)\vee P^{2n(p+1)-1}(p)\vee
P^{4np-2}(p)\vee P^{4np-1}(p).
\]
\end{lemma}
\begin{proof}
The reduced homology $\widetilde{H}_{\ast}(L)$ has three
generators $x,u$ and $v$ in dimensions $2n-1, 2np-2$ and $2np-1$,
respectively. The generators $u$ and $v$ are connected by the
Bockstein homomorphism, that is, ${\beta (v)=u}$. Another way of
looking at the space $L$ is by considering the cofibration
sequence
\begin{equation}
\label{cofibL}
 \ddr{P^{2np-2}(p)}{g}{S^{2n-1}}{}{L}
\end{equation}
that also defines $(2np-1)$-skeleton of $F_2(n)$. Smashing
cofibration sequence \eqref{cofibL} with $L$, there is a
cofibration sequence
\begin{equation}
\label{cofibl2} \ddr{\Sigma P^{2np-2}(p)\wedge L}{\Sigma g\wedge
L}{\Sigma S^{2n-1}\wedge L}{}{\Sigma L\wedge L}.
\end{equation}
Our aim is to show that $\Sigma g\wedge L\simeq *$.

Pinching the bottom sphere of $L$ to the based point, we get the
pinch map
\[
\dr{q\colon L}{}{P^{2np-1}(p)}
\]
that results in the following commutative diagram
\[
\xymatrix{ \Sigma P^{2np-2}(p)\wedge L \ar[d]^{\Sigma
P^{2np-2}(p)\wedge q}\ar[rr]^-{\Sigma g\wedge L} & &\Sigma
S^{2n-1}\wedge L \ar[d]^{\Sigma S^{2n-1}\wedge q}\\
\Sigma P^{2np-2}(p)\wedge P^{2np-1}(p) \ar[rr]^-{\Sigma g\wedge
P^{2np-1}(p)} & & \Sigma S^{2n-1}\wedge P^{2np-1}(p).}
\]
Since $\Sigma^2 g\simeq\Sigma^2 f\simeq *$, it follows that
$\Sigma g\wedge P^{2np-1}(p)\simeq *$. Therefore our problem
reduces to showing that $\gamma \colon \Sigma P^{2np-2}\wedge
L\lra S^{4n-1}$ is null homotopic. Odd sphere $S^{4n-1}$ is an
$H$-space, so it suffices to show that
\[
\Sigma\gamma\colon \Sigma^2 P^{2np-2}\wedge L\lra \Sigma S^{4n-1}
\]
is null homotopic. Suspending cofibration sequence \eqref{cofibl2}
once more, there is a factorisation of the map $\Sigma \gamma$
through the trivial map $\Sigma^2g\wedge L$:
\[
\xymatrix{\Sigma^2 P^{2np-2}\wedge L \ar[r]^-{\Sigma^2g\wedge L}
\ar[dr]^{\Sigma\gamma} & \Sigma^2 S^{2n-1}\wedge L \ar[d]\\
& \Sigma^2 S^{2n-1}\wedge S^{2n-1}.}
\]
Therefore $\Sigma\gamma\simeq *$. This implies that $\gamma$ is
trivial as well.
\end{proof}

To find a left homotopy inverse of
$\dr{\Omega\omega\colon\Omega\Sigma R} {}{\Omega\Sigma L}$ we
require the following two lemmas, proved in \cite{CMNe}.
\begin{lemma}
\label{sigmainverse} A map $\dr{g\colon\Omega X}{}{ Z}$ has a left
homotopy inverse if and only if $\dr{\Sigma g\colon\Sigma\Omega
X}{}{\Sigma Z}$ has a left homotopy inverse.
\end{lemma}
\begin{lemma}
\label{retractP} If $f\colon X\lra Y$ is a map between wedges of
spheres and Moore spaces which induces a mod~$p$ homology
isomorphism, then $f$ has a left homotopy inverse.
\end{lemma}
\begin{lemma}
The homotopy fibration $ \Omega\Sigma R
\stackrel{\Omega\omega}{\longrightarrow}\Omega\Sigma
 L\longrightarrow F$ is a trivial principal fibration, which is to say there exists
 a homotopy decomposition $\Omega\Sigma L\simeq F\times\Omega\Sigma R.$
\end{lemma}
\begin{proof}
To prove the Lemma it suffices to show that there is a left
homotopy inverse of $\Omega\omega$, which by Lemma
\ref{sigmainverse}, exists if $\Sigma\Omega\omega$ has a left
homotopy inverse. By James' theorem, $\Sigma\Omega\Sigma
R\simeq\bigvee^\infty_{i=1}\Sigma R^{(i)}$. Applying the homotopy
equivalences
\[
{S^n\wedge S^m\simeq S^{m+n}}, \  S^{m}\wedge P^n(p^r)\simeq
P^{m+n}(p^r)
 \]
and
\[
P^m(p^r)\wedge P^n(p^r)\simeq P^{m+n}(p^r)\vee P^{m+n-1}(p^r)
\]
to $\bigvee^\infty_{i=1}\Sigma R^{(i)}$, the space
$\Sigma\Omega\Sigma R$ decomposes as a wedge of spheres and Moore
spaces. Similarly, James' theorem gives $\Sigma\Omega\Sigma
L\simeq \bigvee^\infty_{i=1} \Sigma L^{(i)}$. Using
Theorem~\ref{properties}~(iv) and Lemma~\ref{sll}, there is a
homotopy equivalence $\Sigma\Omega\Sigma L\simeq\Sigma
L\bigvee\Bigl( \bigvee_\alpha S^{n_\alpha}\Bigr)\bigvee\Bigl(
\bigvee_\beta P^{n_\beta}(p^r)\Bigr)$. In Proposition
\ref{omegaomega*} we have shown that
${\Omega\omega\colon\dr{\Omega\Sigma R}{}{\Omega \Sigma L}}$
induces a Hopf algebra monomorphism; hence
$(\Sigma\Omega\omega)_*$ does as well. Satisfying the assumption
of Lemma~\ref{retractP}, the map $\Sigma\Omega\omega\colon
\Sigma\Omega\Sigma R\lra \Sigma\Omega\Sigma L$ has a left homotopy
inverse.
\end{proof}
\begin{cor}
\label{homologyF} As a Hopf algebra over the Steenrod algebra,
\[
H_*(F)\cong S(x_{2n-1}, u_{2np-2}, v_{2np-1}).
\]
The action of the Bockstein is given by $\beta (v)=u$, while the
action of other Steenrod powers is trivial.
\end{cor}
\begin{proof}
Consider the principal fibration ${\Omega\Sigma R
\stackrel{\Omega\omega}{\lra}\Omega\Sigma L\lra F}$. Since
$H_*(\Omega\Sigma L)=U\L$ is free $H_*(\Omega\Sigma R)=U[\L,\L
]$-module, the Eilenberg-Moore spectral sequence which converges
to $H_*(F)$ collapses. That is,\newline
$E^2=\mathop{\mathrm{Tor}}^{U[\L,\L]}(\mathbb{Z}/p,U\L)=\mathbb{Z}/p\otimes_{U[\L,\L]}
U\L =U\L_{ab}\langle x,u,v\rangle = S(x,u,v).$ Hence, $H_*(F)\cong
S(x,u,v)$ as a coalgebra.

Use the section $s: F\lra \Omega\Sigma L$ to define the $H$-space
multiplication on $F$ by the composite $F\times F\lra \Omega\Sigma
L\times \Omega\Sigma L\lra \Omega\Sigma L\lra F$. Then
$H_*(F)\cong S(x,u,v)$ as a Hopf algebra.
\end{proof}
Finally, we want to bring $F_2(n)$ into the picture. Since
$F_{2}(n)$ is defined as a homotopy pullback of $H$-spaces and
$H$-maps, it is itself an $H$-space. Thus the Hopf construction
can be applied to $F_{2}(n)$ inducing the following diagram of
homotopy fibration sequences
\begin{equation}
\label{hopfconstruction} \xymatrix{
 \Omega N \ar[r] \ar[d] & \Omega\Sigma L \ar[r]^{\partial} \ar[d]_{\Omega\Sigma i} &
  F_{2}(n)\ar[r] \ar@{=}[d] & N \ar[r] \ar[d] & \Sigma L \ar[d]_{\Sigma i} \\
\Omega(F_{2}(n)\ast F_{2}(n)) \ar[r] & \Omega\Sigma F_{2}(n)
\ar[r] & F_{2}(n) \ar[r] & F_{2}(n)\ast F_{2}(n) \ar[r]^u & \Sigma
F_{2}(n),}
\end{equation}
where $i\colon L\lra F_2(n)$ is the canonical inclusion of the
$(2np-1)$-skeleton, and $N$ is the homotopy pullback of the maps
$u$ and $\Sigma i$.

\begin{lemma}
The map $\Omega\Sigma i\colon\Omega\Sigma L\longrightarrow
\Omega\Sigma F_{2}(n)$ in diagram~\eqref{hopfconstruction} has a
left homotopy inverse $s\colon\Omega\Sigma F_2(n)\lra\Omega\Sigma
L$.
\end{lemma}
\begin{proof}
By Theorem \ref{properties} (v), the space $\Sigma L$ retracts off
$\Sigma F_2(n)$. Therefore, $\Sigma i$ has a left homotopy
inverse, and hence $\Omega\Sigma i$ does as well.
\end{proof}

\begin{pro}
There is a homotopy equivalence
\[
g\colon F_{2}(n)\lra F.
\]
\end{pro}

\begin{proof}
Define a map $g: F_{2}(n)\longrightarrow F$ by the composite
\[
\dddr{F_{2}(n)}{E}{\Omega\Sigma F_{2}(n)}{s}{\Omega\Sigma
L}{\partial}{F}.
\]
In homology, $E_*$ is the canonical inclusion of $H_*(F_2(n))$
into $T(H_*(F_2(n)))$; being a loop map by construction, $s$
induces an algebra map $s_*$; while by the construction
$\partial_*$ is a Hopf algebra monomorphism. Therefore,
$g_*=\partial_*\circ s_*\circ E_*$ is a map of algebras. By
Corollary~\ref{homologyF}, the spaces $F_2(n)$ and $F$ have the
same homology $S(x_{2n-1}, u_{2np-2}, v_{2np-1})$ as Hopf algebras
over the Steenrod algebra. Since in homology $g$ induces an
isomorphism in degrees $2n-1,$ $2np-2$ and $2np-1$, the map $g_*$
is an isomorphism. Thus $g$ is a homotopy equivalence.
\end{proof}

\begin{proof}[Proof of Proposition \ref{fibre}]
Identifying in homotopy fibration sequence~\eqref{fibre2} the
fibre $F$ with the homotopy equivalent space $F_2(n)$, we obtain
the fibration sequence stated by the Proposition.
\end{proof}

\subsection{The homotopy associativity and homotopy commutativity
of $F_{2}(n)$}

To begin the argument proving the homotopy associativity and
homotopy commutativity of $F_{2}(n)$, we recall a pair of
constructions from classical homotopy theory.

Let $X$ be a topological space, and $i_{X}$ the inclusions of $X$
into the wedge ${X\vee X}$. Looping, we can take the Samelson
product $[\Omega i_{X}, \Omega i_X]$. Adjointing gives the
Whitehead product $[\zeta_{X},\zeta_X]$, where
${\zeta_{X}=i_{X}\circ ev}$ and $ev$ is the canonical evaluation
map $ev\colon \Sigma \Omega X\longrightarrow X$. A classical
result in homotopy theory asserts that there is a homotopy
fibration
\begin{equation}
\label{classical} \ddr{\Sigma\Omega X \wedge \Omega
X}{[\zeta_{X},\zeta_X]}{X\vee X}{i}{X\times X}
\end{equation}
where $i$ is the inclusion. The \textit{universal Whitehead
product of} $X$ is defined as the composition
\[
\Sigma\Omega X \wedge \Omega
X\stackrel{[\zeta_{X},\zeta_{X}]}{\lra} X\vee
X\stackrel{\nabla}{\lra} X
\]
where $\nabla$ is the fold map. It has the universal property that
any Whitehead product on $X$ factors through
$[\zeta_{X},\zeta_{X}]$.

The linchpin in showing that $F_{2}(n)$ is \hahc is the following
theorem proved by Theriault \cite{Th2}.
\begin{thm}
\label{htpy assoc&commut} Let $\Omega
B\stackrel{\partial}{\longrightarrow} F \longrightarrow
E\longrightarrow B$ be a homotopy fibration sequence in which
$\partial$ has a right homotopy inverse. Suppose that there is a
homotopy commutative diagram
\[
\xymatrix{\Sigma\Omega B\wedge\Omega B \ar[r] \ar[d] & B\vee B
\ar[d]^{\nabla} \\
E \ar[r] & B, }
\]
where the upper composite in the square is the universal Whitehead
product of $B$. Then the multiplication on $F$ induced by the
retraction off $\Omega B$ is both \hahc and the connecting map
$\partial$ is an $H$-map.
\end{thm}

Coming back to our case, we consider the fibration sequence
\[
\dddr{\Omega\Sigma L}{\partial}{F_{2}(n)}{}{\Sigma
R}{\omega}{\Sigma L}
\]
from Proposition \ref{fibre}. In this case the external universal
Whitehead product is a map $\Psi\colon\Sigma\Omega\Sigma L\wedge
\Omega\Sigma L \longrightarrow \Sigma L\vee \Sigma L$.

\begin{lemma}
\label{p-sum} The external universal Whitehead product
\[
\Psi\colon\Sigma\Omega\Sigma L\wedge \Omega\Sigma L\lra \Sigma
L\vee\Sigma L
\]
is homotopic to a sum of a Whitehead product and mod~$p$ Whitehead
products.
\end{lemma}
\begin{proof}
Using James' theorem Theorem \ref{properties}~$(iv)$ and
Lemma~\ref{sll}, $ \Sigma\Omega\Sigma L\wedge\Omega\Sigma L \simeq
\Sigma M$, where $M$ is a wedge of spheres and Moore spaces. Since
$\Omega\Psi$ has a left homotopy inverse, as fibration
\eqref{classical} splits when looped, it is an inclusion in
homology. Furthermore, the Hurewicz image of each sphere or the
mod~$p$ Hurewicz image of each Moore space summand $P$ of $M$
under the composite
\[
\dddr{\theta\colon M}{E}{\Omega\Sigma
M}{\simeq}{\Omega(\Sigma\Omega\Sigma L\wedge\Omega\Sigma
L)}{\Omega\Psi}{\Omega(\Sigma L\vee\Sigma L)}
\]
is a bracket in $\mathcal{L}\langle V\rangle$, where
$H_{\ast}(\Omega(\Sigma L\vee\Sigma L))\cong H_*(\Omega\Sigma
(L\vee L))\cong U\mathcal{L}\langle V\rangle$ for
${V=\widetilde{H}_*(L\vee L)}$. Using the identity and mod~$p$
Bockstein maps on each summand of $L\vee L$, it is clear that
there exists a mod~$p$ Samelson product on $\Omega(L\vee L)$ which
has the same Hurewicz image as $P$. Summing these mod~$p$ Samelson
products, one for each summand of $M$, gives a map ${\lambda\colon
M\lra\Omega(L\vee L)}$ with the property that
$\lambda_*=\theta_*$. Each mod~$p$ Samelson product factors
through the loop of the universal Samelson product so $\lambda$
lifts to a map $\lambda'\colon M\lra \Omega(\Sigma\Omega L\wedge
L)$ with ${\lambda\simeq \Omega\Psi\circ\lambda'}$. Extend
$\lambda'$ to $\dr{\overline{\lambda}\colon\Omega\Sigma
M}{}{\Omega(\Sigma\Omega P\wedge\Omega P)}$. As
$\lambda_*=\theta_*$, we get
$(\Omega\Psi\circ\overline{\lambda})_*=(\Omega\Psi)_*$. As
$(\Omega\Psi)_*$ is a monomorphism, we must have
$(\overline{\lambda})_*$ is an isomorphism. Hence
$\overline{\lambda}$ is a homotopy equivalence. Taking adjoints
then proves the Lemma.
\end{proof}
\begin{lemma}
\label{lift} There is a lift
\[
\xymatrix{ & \Sigma\Omega\Sigma L\wedge\Omega\Sigma L
\ar@{-->}[dl]\ar[d]^{\nabla\circ\Psi}\\
\Sigma R \ar[r]^{\omega} & \Sigma L}
\]
of the universal Whitehead product of $\Sigma L$ to $\Sigma R$.
\end{lemma}
\begin{proof}
Lemma \ref{p-sum} shows that the universal Whitehead product on
$\Sigma L$ is homotopic to a sum of Whitehead products and mod~$p$
Whitehead products. The Whitehead product defines the Lie bracket
on $\Sigma L$. The set $W$ from Proposition~\ref{basis} consists
of Whitehead products and mod~$p$ Whitehead products which form a
Lie basis for $[\mathcal{L},\mathcal{L}]$. So mod~$p$ Whitehead
products from the universal Whitehead product can be rewritten as
a linear combination of basis elements.
\end{proof}
\begin{pro}
\label{Fac} $F_{2}(n)$ is homotopy associative and homotopy
commutative.
\end{pro}
\begin{proof}
Applying Theorem \ref{htpy assoc&commut} to the fibration sequence
\[
\dddr{\Omega\Sigma L}{\partial}{F_{2}(n)}{*}{\Sigma
R}{\omega}{\Sigma L}
\]
and using Lemma \ref{lift}, the Proposition follows.
\end{proof}

\subsection{A universal property of $F_{2}(n)$}

In this section we show that $F_2(n)$ satisfies the universal
property in the category of homotopy associative, homotopy
commutative $H$-spaces. Let $f\colon L\lra Z$ be a map into a
homotopy associative, homotopy commutative $H$-space. We want to
show that there is a unique multiplicative extension
$\overline{f}\colon F_2(n)\lra Z$ of $f$.

\begin{pro}
\label{hcom} Let $h\colon\Omega\Sigma L \longrightarrow Z$ be an
$H$-map into a homotopy commutative and homotopy associative
$H$-space $Z$. Then it factors through $\partial\colon\Omega\Sigma
L\lra F_2(n)$.
\end{pro}

\begin{proof}
Denote by $g$ a right homotopy inverse of the $H$-map $\partial$
and by $e$ the homotopy equivalence ${e\colon \dr{F_2(n)}\times
\Omega\Sigma R}{g\cdot\Omega\omega}{\Omega\Sigma L}$. Define two
maps $a,b\colon\Omega\Sigma L\lra \Omega\Sigma L$ as the
composites
\[
\dddr{a\colon\Omega\Sigma L}{e^{-1}}{F_2(n)\times \Omega\Sigma
R}{\pi_1}{F_2(n)}{g}{\Omega\Sigma L}
\]
\[
\dddr{b\colon\Omega\Sigma L}{e^{-1}}{F_2(n)\times \Omega\Sigma
R}{\pi_2}{\Omega\Sigma R}{\Omega\omega}{\Omega\Sigma L}.
\]
In the following diagram
\[
\xymatrix{\Omega\Sigma L\ar[r]^-\Delta \ar[rd]^{e^{-1}}&
\Omega\Sigma L\times \Omega\Sigma L\ar[r]^-{(a,b)} & \Omega\Sigma
L
\times\Omega\Sigma L \ar[r]^-{\mu} & \Omega\Sigma L\\
& F_2(n)\times \Omega\Sigma R \ar[r]^-{(g,\Omega\omega)} &
\Omega\Sigma L\times\Omega\Sigma L\ar[u] &}
\]
the top row composition is $a+b$, while the bottom row is the
identity on $\Omega\Sigma L$. The commutativity of the diagram
gives $\mathrm{Id}_{\Omega\Sigma L}\simeq a+b$.

Being an $H$--map, $h$ is determined by its restrictions on the
both of the factors of $\Omega\Sigma L$, that is, $h\simeq h\circ
\mathrm{Id}_{\Omega\Sigma L}\simeq h\circ(a+b)\simeq h\circ a
+h\circ b$.

Showing that the composite
\[
\ddr{\Omega\Sigma R}{\Omega\omega}{\Omega\Sigma L}{h}{Z}
\]
is null homotopic, we have $h\circ b\simeq *$ and hence $h\simeq
h\circ a$ that proves the Proposition.

As a composite of $H$-maps, $h\circ\Omega\omega$ is itself an
$H$-map. Therefore by the James construction, it is uniquely
determined by its restriction to $R$. The composite
\[
\ddr{R}{E}{\Omega\Sigma R}{\Omega\omega}{\Omega\Sigma L}
\]
is a wedge sum of Samelson products on $R$ as it is the adjoint of
a wedge sum of Whitehead products ${\omega\colon\dr{\Sigma
R}{}{\Sigma L}}$. Being an $H$-map, $h$ preserves Samelson
products. Therefore the wedge of Samelson products
$\Omega\omega\circ E$ on $R$ composed with $h$ into the homotopy
commutative $H$-space $Z$ is trivial.
\end{proof}
Consider the fibration sequence
\[
\dddr{\Omega\Sigma L}{\partial}{F_2(n)}{*}{\Sigma
R}{\omega}{\Sigma L}.
\]
Define a map $\overline{f}$ as the composite
\[
\overline{f}\colon
F_{2}(n)\stackrel{g}{\longrightarrow}\Omega\Sigma
L\stackrel{\tilde{f}}{\longrightarrow} Z,
\]
where $g\colon F_{2}(n)\longrightarrow \Omega\Sigma L$ is a right
homotopy inverse of the map ${\partial\colon\Omega\Sigma
L\longrightarrow F_{2}(n)}$, and $\tilde{f}\colon\Omega\Sigma
L\longrightarrow Z$ is the canonical extension of $f$ given by the
James construction.

\begin{thm}
\label{universalityf2} Let $Z$ be a \hahc $H$-space, and $f\colon
L\longrightarrow Z$ a given map. Then there exists an extension to
a unique up to homotopy $H$-map ${\overline{f}\colon
F_{2}(n)\longrightarrow Z}$.
\end{thm}
\begin{proof}
Let us consider the map $\overline{f}$ given by
\[
\overline{f}\colon
F_{2}(n)\stackrel{g}{\longrightarrow}\Omega\Sigma
L\stackrel{\tilde{f}}{\longrightarrow} Z,
\]
as a candidate for the multiplicative extension of $f$. Notice
that there is a commutative diagram
\begin{equation}
\label{extension}
\xymatrix{L \ar[d]^E\ar[dr]^f &\\
\Omega\Sigma L \ar[r]^{\widetilde{f}}\ar[d]^{\partial} &
Z\\
F_2(n) \ar[ur]^{\overline{f}}}
\end{equation}
since the upper triangle commutes by the James construction, while
the commutativity of the lower diagram is given by
Proposition~\ref{hcom}. Diagram~\eqref{extension} ensures that
$\overline{f}$ is an extension of $f$.

Now, we shall show that $\overline{f}$ is an $H$-map by showing
that the diagram
\begin{equation}
\label{Hmap} \xymatrix{F_2(n)\times F_2(n) \ar[r]^-{g\times g}
\ar[dr]^-{\tilde{f}g\times\tilde{f}g} & \Omega\Sigma
L\times\Omega\Sigma L
\ar[r]^-{\mu}\ar[d]^-{\tilde{f}\times\tilde{f}} & \Omega\Sigma
L\ar[r]^-{\partial}\ar[d]^-{\tilde{f}} &
F_2(n)\ar[d]^-{\tilde{f}g}\\
& Z\times Z \ar[r]^-{\mu} & Z\ar@{=}[r] & Z }
\end{equation}
commutes. The left triangle commutes by definition; the middle
square commutes since $\tilde{f}$ is an $H$-map; and the
commutativity of the right square is given by
Proposition~\ref{hcom}. Summing this up, diagram~\eqref{Hmap}
commutes.

Finally we are left to show that ${\overline{f}=\tilde{f}\circ
g\colon F_{2}(n)\longrightarrow Z}$ is the unique $H$-map
extending $f\colon L\longrightarrow Z$. To prove that we use the
uniqueness of $\tilde{f}$ asserted by the James construction and
the result of Theorem \ref{htpy assoc&commut} which shows that the
connecting map $\partial\colon\dr{\Omega\Sigma L}{}{F_2(n)}$ is an
$H$-map. Let $\overline{f}, \overline{g}$ be two extensions
$\dr{F_2(n)}{\overline{f}, \overline{g}}{Z}$ of the map
$f\colon\dr {L}{}{Z}$ which are $H$-maps. Precompose both maps
with the $H$-map $\partial\colon\dr{\Omega\Sigma L}{}{F_2(n)}$. We
obtain two multiplicative extensions $\ddr{\Omega\Sigma
L}{\partial}{F_2(n)}{\overline{f},\ \overline{l}}{Z}$ of
$f\colon\dr{L}{}{Z}$. By the uniqueness of an $H$-map
$\Omega\Sigma L\lra Z$ extending $f$ in the James construction, it
follows that $\overline{f}\circ\partial
\simeq\overline{l}\circ\partial$. Precomposing both compositions
with the right homotopy inverse $g\colon\dr{F_2(n)}{}{\Omega\Sigma
L}$ of the map $\partial$, we get
\[
\overline{f}\circ\partial\circ g\simeq
\overline{l}\circ\partial\circ g.
\]
Hence
\[
\overline{f}\simeq \overline{l}
\]
and the uniqueness assertion is proved. This finishes the proof of
the Theorem.
\end{proof}

The following theorem summarises Proposition \ref{Fac} and Theorem
\ref{universalityf2}.
\begin{thm}
\label{unispLa} $F_2(n)$ is a universal space of $L$.
\end{thm}

The uniqueness assertion of the theorem is powerful. It can be
used to show that two $H$-maps from $F_{2}(n)$ to a \hahc
$H$-space are homotopic by comparing their restriction to $L$.
This property is the foundation stone in the application we
consider in the next section.

\section{The $d_{1}$-differential in the $EHP$ spectral sequence}
\subsection{Extension of a formula for the $d_1$-differential to $F_2(np)$}

The problem we want to discuss in this section is that of
calculating the unstable homotopy groups of spheres using the
$EHP$ spectral sequence and applying the universality of $F_2(n)$.

Recall from the Introduction that the $d_1$-differential in the
$EHP$ spectral sequence is induced by the composition
$\ddr{\Omega^3 S^{2np+1}}{\Omega P}{\Omega
J_{p-1}(S^{2n})}{H}{\Omega S^{2np-1}}$. In~\cite{Gr1} Gray
constructed a map $\varphi_n\colon \Omega^2 S^{2np+1}
\longrightarrow S^{2np-1}$ with the property that $\Omega
\varphi_n= H\circ\Omega P$. Therefore, the $d_1$-differential can
be considered as $d_{1}(x)=\varphi_{n}\cdot x$. Using the
existence of the map $\varphi_n$ and the universality of $\Omega
J_{p-1}(S^{2n})$ when $p\geq 3$, Gray \cite{Gr2} developed a
formula for $d_1$ applicable to elements which are not double
suspensions.

One would like to have a formula for the $d_{1}$-differential
applicable to its whole domain. Since $F_{2}(n)$ is a third
approximation of the double suspension, after $S^{2n-1}$ and
$\Omega J_{p-1}(S^{2n})$, we would like to extend the description
of $d_{1}$ from $\Omega J_{p-1}(S^{2np})$ to $F_{2}(np)$.

There are a few points which should be commented upon. Recall that
$W_n$ denotes the homotopy fibre of the double suspension
$E^2\colon S^{2n-1}\lra \Omega^2 S^{2n+1}$. It has been shown in
\cite{Gr3} that $W_n$ is a loop space, namely that there exists a
fibration sequence
\[
W_{n}\longrightarrow
S^{2n-1}\stackrel{E^{2}}{\longrightarrow}\Omega^{2}S^{2n+1}\longrightarrow
BW_{n},
\]
where $BW_n$ is defined as the homotopy fibre of $\varphi_n$,
\[
BW_{n}\longrightarrow
\Omega^{2}S^{2np+1}\stackrel{\varphi_{n}}{\longrightarrow}
S^{2np-1}.
\]

The following lemma states the existence of two fibration
sequences involving $F_2(n)$ which are analogous to the two
classical $EHP$ fibrations.
\begin{lemma}
There exist homotopy fibration sequences
\[
 \dddr{W_{np}}{P}{F_2(n)}{E}{\Omega^2
S^{2n+1}}{H}{BW_{np}}
\]
and
\[
\dddr{\Omega S^{2np-1}\{
p\}}{P}{S^{2n-1}}{E}{F_2(n)}{H}{S^{2np-1}\{p\}.}
\]
\end{lemma}
\begin{proof}
The homotopy pullback defining $F_2(n)$ can be extended to the
homotopy pullback diagram
\[
\xymatrix{F_2(n) \ar[r]\ar[d] & S^{2np-1} \ar[r]\ar[d]^{E^2} &
J_{p-1}(S^{2n})\ar@{=}[d] \\
\Omega^2 S^{2n+1}\ar[r]^{\Omega H} \ar[d] & \Omega^2 S^{2np+1}
\ar[r]^P\ar[d] & J_{p-1}(S^{2n})\\
BW_{np} \ar@{=}[r] & BW_{np}. & }
\]
of homotopy fibration sequences. The left column of the diagram is
the first fibration sequence whose existence is claimed in the
Lemma.

The second homotopy fibration sequence of the Lemma is obtained
from the following homotopy pullback
\[
\xymatrix{S^{2n-1}\ar[r]^-E \ar@{=}[d] & \Omega J_{p-1}(S^{2n})
\ar[r]^-H\ar[d] & \Omega S^{2np-1}\ar[d]\\
S^{2n-1}\ar[r]^-E & F_2(n)\ar[r]^-H & S^{2np-1}\{p\},}
\]
whose existence is proved by Selick in \cite[Theorem 10]{Se}.
\end{proof}
In our programme for finding a formula for the $d_1$-differential
we consider the map
\[
\alpha\colon \dddr{F_2(np)}{E}{\Omega^2
S^{2np+1}}{\varphi_n}{S^{2np-1}}{E}{F_2(np).}
\]
Recall that $L_{np}$ is the $(2np^2-1)$-skeleton of $F_2(np)$. The
space $F_2(np)$ is an $H$-space, so it has a $p$-power map.
\begin{pro}
\label{alpha}
\begin{enumerate}\renewcommand{\labelenumi}{(\roman{enumi})}
 \item The restriction of
 \[
 \alpha\colon \dddr{F_2(np)}{E}{\Omega^2
S^{2np+1}}{\varphi_n}{S^{2np-1}}{E}{F_2(np)}
 \]
to $L_{np}$ is homotopic to the restriction of the $p$-power map
to $L_{np}$.
 \item If $\alpha$ is an $H$-map, then $\alpha\simeq p$.
\end{enumerate}
\end{pro}
\begin{proof}
Look at the difference $D$ of $\alpha$ and the $p$-power map
restricted to the $(2np^2-1)$-skeleton of $F_2(np)$. Since
$(F_2(np))_{(2np^2-1)}$ is homotopy equivalent to the
$(2np^2-1)$-skeleton of $\Omega^2 S^{2np+1}$ the map $\alpha$
restricted to $(F_2(np))_{(2np^2-1)}$ is homotopic to the
composite
\[
\widetilde{\alpha}\colon
\ddr{(\Omega^2S^{2np+1})_{(2np^2-1)}}{\varphi_n}{S^{2np-1}}{E}{F_2(np)}.
\]
Theriault has proved \cite{Th3} that $E^2\circ\varphi_n\simeq p$.
Hence the difference between $E\circ\alpha$ and the $p$ power map
is null homotopic. Therefore $D$ factors through $W_{np^2}$, the
homotopy fibre of $E\colon F_2(np)\lra \Omega^2 S^{2np+1}$.
Applying the homology Serre spectral sequence to the fibration
sequence $\ddr{W_{np^2}}{}{F_2(np)}{E}{\Omega^2S^{2np+1}}$, we
have $W_{np^2}$ is $2np^3-4$ connected, and therefore $D$ is null
homotopic.

The second assertion of the Proposition follows from part (i)
since $F_2(np)$ is the universal space of $L_{np}$. Namely, all
relevant maps are $H$-maps and their composite is uniquely
determined by the restriction to $L_{np}$.
\end{proof}
 The following Theorem gives the main result, which is obtained by applying the
 universality of $F_2(n)$ to the $EHP$ spectral sequence
 calculation of the unstable homotopy groups of spheres.

\begin{thm}
\label{formulad{1}2} The composite
\[
F_2(np)\stackrel{E}{\lra}\Omega^2S^{2np+1}\stackrel{\varphi_n}
{\lra}S^{2np-1}\stackrel{E}{\lra} F_2(np)
\]
is the $p$-power map if either
\begin{enumerate}
 \item [(i)] $x\in\pi_*(F_2(np))$ is an element which is a lift through ${H\colon
 F_2(np)\lra S^{2np^2-1}\{p\}}$ of an element in the image of ${E\colon P^{2np^2-2}(p)\lra \Omega
 S^{2np^2-1}\{p\}}$
\end{enumerate}
or
\begin{enumerate}
 \item[(ii)] the map $\alpha$ is an $H$-map.
\end{enumerate}
\end{thm}
\begin{proof}
First we show that if $H(x)=E(u)$, then $x\in \Im
(\pi_*(L_{np})\longrightarrow \pi_*(F_2(np))$. The assumptions can
be presented in the form of the diagram
\[
\xymatrix{\Omega F_2(np) \ar[r]^-{\Omega H} & \Omega
S^{2np^2-1}\{p\} \ar[r]^-P & S^{2np-1}\\
S^m \ar[u]_x \ar[r]^-u & P^{2np^2-2}(p)\ar[u]_E.  &}
\]
Let $\omega_{np}$ be the composite $\ddr{P^{2np^2-2}(p)}{E}{\Omega
S^{2np^2-1}\{ p\} }{P}{S^{2np-1}}$. Then
\[
\omega_{np}\circ u\simeq 0
\]
since it factors through two consecutive terms of the homotopy
fibration sequence
\[
\ddr{\Omega F_{2}(np)}{\Omega H}{\Omega S^{2np^2-1}\{
p\}}{P}{S^{2np-1}}.
\]
Therefore $u$ factors through the fibre $F_{\omega_{np}}$ of
$\omega_{np}$. Putting this all together we obtain the following
diagram
\[
\xymatrix{ & F_{\omega_{np}} \ar[rr]\ar[d] & & \Omega F_2(np)\ar[d]^{\Omega H}\\
S^m \ar[r]^-u \ar@{-->}[ur]^{\widetilde{x}} & P^{2np^2-2}(p)
\ar[rr]^-{E} \ar[d]^{\omega_{np}} & & \Omega
S^{2np^2-1}\{ p\} \ar[d]^P\\
& S^{2np-1}\ar@{=}[rr]& & S^{2np-1},}
\]
where $\widetilde{x}$ is a lift of $u$ to $F_{\omega_{np}}$.

It follows from the defining property of a pullback that the map
$x$ factors as
\[
x\colon \ddr{S^m}{\widetilde{x}}{F_{\omega_{np}}}{}{\Omega
F_2(np)}.
\]
We claim that the map $\dr{F_{\omega_{np}}}{}{\Omega F_2(np)}$
factors through $\Omega L_{np}$. To see this, let $X$ be the fibre
of the pinch map $q\colon\dr{L_{np}}{}{P^{2np^2-1}(p)}$ to the top
Moore space. Then there is a fibration sequence
\[
\Omega P^{2np^2-1}(p) \lra X \lra L_{np} \lra P^{2np^2-1}(p).
\]
Starting with the diagram
\[
\xymatrix{L_{np}\ar[r]\ar[d] & F_2(np)\ar[d]\\
P^{2np^2-1}(p)\ar[r] & S^{2np^2-1}\{ p\}}
\]
that includes the top Moore space of $L_{np}$ into $S^{2np^2-1}\{
p\}$, there is an induced map of homotopy fibrations
\[
\xymatrix{\Omega P^{2np^2-1}(p)\ar[r]\ar[d] & \Omega S^{2np^2-1}\{ p\}\ar[d]\\
X\ar[r]\ar[d] & S^{2np-1}\ar[d]\\
L_{np}\ar[r] & F_2(np).}
\]
Including the bottom cell $S^{2np-1}$ into $X$, let us consider
the composition
\[
\dddr{f\colon
P^{2np^2-2}(p)}{\omega_{np}}{S^{2np-1}}{i}{X}{}{L_{np}}.
\]
By looking at the map
\[
\xymatrix{P^{2np^2-2}(p)\ar[d]^E
\ar[r]^-{\text{attaching}}_-{\text{map}} & S^{2np-1}
\ar[r]\ar[d]^{=} & L_{np}\ar[r]\ar[d] & P^{2np^2-1}(p)
\ar[d]\\
 \Omega S^{2np^2-1}\{ p\}\ar[r]^-P & S^{2np-1} \ar[r] & F_2(np) \ar[r] & S^{2np^2-1}\{ p\}}
\]
between a cofibration homotopy sequence and a fibration sequence,
we conclude that the map $\omega_{np}$ is homotopic to the
attaching map defining $L_{np}$. Thus the map $f$ is null
homotopic since it factors through the cofibration sequence
\[
P^{2np^2-2}(p)\lra S^{2np-1}\lra L_{np}.
\]
Therefore, there exists a lift $E\colon P^{2np^2-2}(p)\lra \Omega
P^{2np^2-1}(p)$ closing the diagram
\[
\xymatrix{ & & \Omega P^{2np^2-1}(p)\ar[d]\\
P^{2np^2-2}(p)\ar@{-->}[rru]^E \ar[r]^{\omega_{np}} & S^{2np-1}
\ar@{^{(}->}[r]^i & X.}
\]

This lift allows us to take the juxtaposition of the two fibration
diagrams
\begin{equation}
\label{factorisationdiagram} \xymatrix{F_{\omega_{np}}
\ar[r]\ar[d] & \Omega L_{np}
\ar[r]\ar[d] & \Omega F_2(np) \ar[d]\\
P^{2np^2-2}(p) \ar[r]^-E\ar[d]^{\omega_{np}} & \Omega
P^{2np^2-1}(p) \ar[r]\ar[d] &
\Omega S^{2np^2-1}\{ p\}\ar[d]\\
S^{2np-1}\ar@{^{(}->}[r] & X \ar[r] & S^{2np-1}}
\end{equation}
Diagram \eqref{factorisationdiagram} shows that there is a
factorisation of the map $\dr{F_{\omega_{np}}}{}{\Omega F_2(np)}$
through $\Omega L_{np}$. Now both statements of the Theorem follow
from Proposition~\ref{alpha}.
\end{proof}
\begin{rem}
Passing to the homotopy groups under the assumptions of
Theorem~\ref{formulad{1}2} the formula for the $d_1$-differential
takes the form
\begin{equation}
\label{formularem1} Ed_1(Ex)=p\cdot x.
\end{equation}
\end{rem}

\subsection{Partial results}

\begin{pro}
Restricted to the $(2np^3-4)$-skeleton the composite
\[
F_2(np)\lra\Omega^2 S^{2np+1}\stackrel{\varphi_n}{\lra}
S^{2np-1}\stackrel{E}{\lra}F_2(np)
\]
is homotopic to the $p$-power map restricted to
$(F_2(np))_{(2np^3-4)}$.
\end{pro}
\begin{proof}
First let us notice that $(F_2(np))_{(2np^3-3)}$ is homotopy
equivalent to the $(2np^3-3)$-skeleton of $\Omega^2 S^{2np+1}$.
Therefore, the map
\[
{\overline{\alpha}\colon\ddr{(\Omega^2
S^{2n+1})_{(2np^3-4)}}{\varphi_n}{S^{2np-1}}{E}{F_2(np)}}
\]
is homotopic to the restriction of
\[
\alpha \colon\dddr{F_2(np)}{E}{\Omega^2
S^{2np+1}}{\varphi_n}{S^{2np-1}}{E}{F_2(np)}
\]
to the $(2np^3-3)$-skeleton of $F_2(np)$. Theriault has proved
\cite{Th3} that ${E^2\circ\varphi_n\simeq p}$. Hence the
difference $D$ between $E\circ\alpha$ and the $p$-power map is
null homotopic. Therefore $D$ factors through $W_{np^2}$, the
homotopy fibre of $E\colon F_2(np)\lra \Omega^2 S^{2np+1}$. Since
$W_{np^2}$ is $2np^3-4$ connected, $D$ is null homotopic.
\end{proof}
\begin{cor}
Formula \eqref{formularem1} for the $d_1$-differential in the
$EHP$-spectral sequence is correct up to the
${(2np^3-4)\text{-skeleton}}$ of $\Omega^2 S^{2np+1}$.
\end{cor}
\begin{proof}
It follows immediately from Theorem \ref{formulad{1}2} (ii) since
restricted to the given skeleton of $\Omega^2 S^{2np+1}$ the map
$\alpha$ is an $H$-map.
\end{proof}

\bibliographystyle{amsplain}

\end{document}